\magnification=1200

\tolerance=500

\font \twelvebf=cmbx12

\def \Ext {{\rm Ext\,}}
\def \cExt {{{\cal E}\it xt\,}}
\def\mapright#1{\smash{\mathop{\longrightarrow}\limits^{#1}}}
\def \titre#1{\medbreak \noindent {\bf #1.}\medbreak}
\outer \def \th #1. #2\par{ \medbreak 
\noindent {\bf#1. \enspace} {\sl#2 }\par
\ifdim \lastskip< \medskipamount \removelastskip \penalty55 \medskip \fi}
\def \lign {\hfil \break }
\def\pf{\noindent {\bf Proof.}\enspace \nobreak }
\def \remk#1 {\medbreak \noindent {\bf Remark #1.}\enspace }
\def \remks#1 {\medbreak \noindent {\bf Remarks #1.}\enspace }
\def \exp#1 {\medbreak \noindent {\bf Example #1.}\enspace }
\def \tas {{\hbox {\bf .}}}

\def \bP {{\bf P}}
\def \cB {{\cal B}}
\def \cE {{\cal E}}
\def \cL {{\cal L}}
\def \cO {{\cal O}}
\def \cP {{\cal P}}
\def \cF {{\cal F}}
\def \cJ {{\cal J}}
\def \cM {{\cal M}}
\def \cN {{\cal N}}

\def \bZ {{\bf Z}}
\def \bN {{\bf N}}

\centerline {\twelvebf From coherent sheaves to curves in ${\bf P}^3$}

\vskip 2 cm

\titre {0. Introduction}
Let $k$ be an algebraically closed field, and let 
${\bf P}={\bf P}^3_k$ be the projective 3-space over $k$. We
denote by $R$ the associated polynomial ring  $R=k[X,Y,Z,T]$.

The aim of this paper is to precise the correspondences between equivalence
classes of some objects related to ${\bf P}$ : {\it vector bundles} (i.e.
locally free coherent sheaves), {\it coherent sheaves}, {\it finite
length graded
$R$-modules} and {\it locally Cohen-Macaulay curves} (i.e. closed subscheme 
of pure dimension one with no embedded points). Some of these
correspondences are well-known, others are new.
\vskip 0.3cm
The notion of {\it  pseudo-isomorphism} was introduced in [HMDP] :

\th {D\'efinition 0.1}. Let $A$ be a ring, let $\cN$ and  $\cN'$  be two 
coherent sheaves on
 ${\bf P}^3_A$, flat over $A$  and let 
$f$  be a morphism from 
$\cN$ to $\cN'$.  We say that  $f$ is a  pseudo-isomorphism    if it
induces~:{\hfil \break }  
i) an isomorphism of functors (from the
category of finite type $A$-modules to itself) $H^0(\cN(n)\otimes_A
\tas) \to H^0(\cN'(n)\otimes_A
\tas)$ for all $n \ll 0$,  {\hfil \break }
ii) an isomorphism of functors (from the
category of finite type $A$-modules to the category of the
graded $A[X,Y,Z,T]$-modules)
$H^1_*(\cN\otimes_A \tas) \to H^1_*(\cN'\otimes_A \tas)$, \lign
 iii) a monomorphism of functors $H^2_*(\cN\otimes_A \tas) \to
H^2_*(\cN'\otimes_A
\tas)$. \lign  
Two  coherent sheaves on
 ${\bf P}^3_A$, flat over $A$, are pseudo-isomorphic if there exists a chain
of 
 pseudo-isomorphisms connecting them  :
$$ \cN= \cN_0 \to \cN_1 \leftarrow \cN_2 \to \cN_3 \leftarrow \cdots \to
\cN_{2p-1}
\leftarrow
\cN_{2p} =
\cN'.$$

This defines an equivalence relation called psi-equivalence. It turns out
that when $A$ is a field the psi-equivalence is an
extension to the set  of  coherent sheaves of (local)
projective dimension $\leq 1$ of the stable equivalence for vector bundles
with $H^2_*=0$ (cf. 2.4).

Now,  using the  bijective correspondence between the set of isomorphism
classes of finite length graded
$R$-modules and the set of stable equivalence classes of vector bundles
with $H^2_*=0$ ([Ho], see also 2.3), and a well-known theorem of Rao ([R]),
we get a new correspondence :

\th {Corollary 2.8}.  The application which maps a curve $C$ to its sheaf
of ideals ${\cal J}_C$ induces a bijective correspondence between the  set 
of biliaison classes of  curves and
the  set  of  coherent sheaves of (local)
projective dimension $\leq 1$  up to  psi-equivalence and up to a twist.

In the third paragraph, we show that there exist rank
2 reflexive
 sheaves in every class of psi-equivalence, and among them,  {\bf minimal} 
rank 2 reflexive sheaves, in the sense that their first  Chern class is
minimal. Moreover Buraggina (cf. [B]) points out that if there are rank 2
vector bundles in the class, they are minimal elements and all the minimal
elements are vector bundles.

\vskip 0.3cm

Then two very natural questions arise, concerning the possible relationships
between the minimal curves in a biliaison class and the minimal rank 2
reflexive sheaves in a class of psi-equivalence (recall that the curves
obtained as the zero locus of a section of a twist of a  rank 2 vector
bundle are exactly the  subcanonical curves (cf. H])):

\th{Question I}. Is a minimal curve of a biliaison class the zero locus of a
section of a twist of a (minimal in the corresponding class of
psi-equivalence) rank 2 reflexive sheaf ?

\th{Question II}. If  a biliaison class contains subcanonical curves, is a 
minimal curve subcanonical ?

It is known since [B] that the answer to Question I is NO (we give in 3.9 an
explicit example of a minimal curve which is not the zero locus of a
section of a twist of a  rank 2 reflexive sheaf) and we have recently proved
that the answer to Question II is YES (cf. [MD]).

 \titre {1. The main equivalence relations}

In this section we will discuss some more or less well-known equivalence
relations.

We say that a vector bundle which is isomorphic to a
finite direct sum
$\bigoplus {\cal O}_{\bf P}(-n_i)$  is  {\it dissoci\'e}.

If ${\cal N}$ is a sheaf on ${\bf P}^3_k$, we denote by 
$H_*^i({\cal N})$ the graded
$R$-module $\bigoplus_{n \in
{\bZ}} H^i({\bf P},{\cal N}(n))$, and by $h^i(\cN)$ the dimension of the
vector space $ H^i({\bf P},{\cal N}(n))$.

\th {Definition 1.1}. Two vector bundles $\cF$ and $\cF'$ on ${\bf P}^3_k$
are stably isomorphic if there  exist two vector bundles
 $\cL$ et $\cL'$ which are dissoci\'e and an 
 isomorphism $\cF\oplus \cL\simeq \cF'\oplus \cL'$.

We will denote by ${\cal
S}tab$ the  set of vector bundles $\cF$ such that
$H_*^2({\cal F})=0$ up to stable equivalence.

\vskip 0.3cm
The notion of {\it  pseudo-isomorphism} was introduced in [HMDP] :

\th {Definition 1.2}. Let $\cN$ and  $\cN'$  two  coherent sheaves on
 ${\bf P}^3_k$  and let 
$f$  be a morphism from 
$\cN$ to $\cN'$. We say that  $f$ is a  pseudo-isomorphism   (shortly a 
{\it psi}) if it induces~:{\hfil \break } 
i) an isomorphism 
$H^0\cN(n) \simeq H^0\cN'(n)$ for all $n \ll 0$,  {\hfil \break }
ii) an isomorphism  $H^1_*\cN \simeq H^1_*\cN'$,{\hfil \break } 
iii) an injective homomorphism $H^2_*\cN \to H^2_*\cN'$.{\hfil \break } 
Two  coherent sheaves $\cN$ and  $\cN'$ are pseudo-isomorphic if there exists
a chain of  {\it psi} connecting them  :
$$ \cN= \cN_0 \to \cN_1 \leftarrow \cN_2 \to \cN_3 \leftarrow \cdots \to
\cN_{2p-1}
\leftarrow
\cN_{2p} =
\cN'.$$

This defines an equivalence relation called psi-equivalence. (In fact, this
definition was given  in a more general context, i.e. for families of
coherent sheaves in the case when $k$ is not a field but a ring.)

\exp {1.3}  Let $\cN$ be a reflexive sheaf of rank $r+1$, let $\cP$ be a
dissoci\'e sheaf of rank $r$ and let
$C$ be a curve such that there exists an exact sequence :
 $0\to \cP \to \cN \to
\cJ_C\to 0$, where $\cJ_C$ is the ideal sheaf of $C$. Then the homomorphism
$\cN \to\cJ_C$ arising from the exact sequence is a {\it psi} of coherent
sheaves of  projective dimension $\leq 1$.

\vskip 0.3cm

The following definition is well-known :

\th {Definition 1.4}. Two (locally Cohen-Macaulay) curves $C$ and $C'$ are
linked if there exists a complete intersection curve $X$ containing $C$ and
$C'$ and such that $ {\cal H}om_{{\cal O}_{\bf P}} ({\cal O}_C, {\cal
O}_X)\simeq {\cal I}_{C',X}$
and
$ {\cal H}om_{{\cal O}_{\bf P}} ({\cal O}_{C'}, {\cal O}_X)\simeq {\cal
I}_{C,X}$.
Two curves $C$ and $C'$ are in the same biliaison class if
they can be connected by a chain of an even number of linkages.

We will denote by ${\cB}il$ the  set of biliaison classes of locally
Cohen-Macaulay curves. 

We denote by ${\cM}_f$ the  set of finite length graded
$R$-modules, modulo isomorphism.

\titre {2. Correspondences between quotient sets}

We will begin by proving that the psi-equivalence retains some kind of
regularity of coherent sheaves.

\th {Lemma 2.1}. Let $f:\cF \to \cF'$ be a psi of coherent sheaves. Then
there exist two dissoci\'e sheaves $\cL$ et $\cL'$  and an exact sequence $0
\to
\cL' \to \cL \oplus \cF \to
\cF' \to 0$. Moreover if
$\cF'$ is locally free and satisfies $H_*^2({\cal F}')=0$, then $\cF$ is
also locally free and $\cF$ and 
$\cF'$ are stably isomorphic.

\pf  There exists $n_0\in \bN$ such that $H^0\cF(n) \simeq H^0\cF'(n)$ is
an isomorphism for all $n < n_0$. Let $L$ be a free
$R$-module of finite type, such that there exists a surjective homomorphism
: $L\to \oplus_{n\geq n_0} H^0\cF'(n)$. Denoting by
$\cL$ the dissoci\'e sheaf associated to $L$ and combining the associated
surjective  homomorphism of sheaves : $\cL\to \cN'$ with
$f$, we get a surjective
homomorphism
: $\cL \oplus \cF
\to
\cF'$ whose kernel is  denoted by $\cL'$. So we have an exact
sequence : $0 \to
\cL' \to \cL \oplus \cF \to
\cF' \to 0$. Taking the associated cohomology exact sequence we get :
$$L\oplus H^0_*\cF\to H^0_*\cF' \to H^1_*\cL' \to H^1_*\cF \to H^1_*\cF' \to
H^2_*\cL'
\to H^2_*\cF \to H^2_*\cF'$$
Since $f$ is a psi and $L\to
H^0_*\cF'$ is surjective, we have $H^1_*\cL'=H^2_*\cL'=0$. Thus it
follows by a theorem of Horrocks [Ho] that $\cL'$ is dissoci\'e.

Now  by local duality, $H_*^2({\cal F}')=0$ implies that  $\Ext^1
(\cF',\cL')=0$ so that  the exact sequence splits and we are done.

\th {Corollary 2.2}. Let $\cF$ and  $\cF'$ two coherent sheaves which are
pseudo-isomorphic. Then the sheaves $\cExt^i_{\cO_\bP}
(\cF,\cO_\bP)$ and $\cExt^i_{\cO_\bP}(\cF',\cO_\bP)$ are isomorphic for
$i\geq 2$.

\pf Suppose that we have a psi $f:\cF \to \cF'$. Then by 2.1, there exist two
dissoci\'e sheaves $\cL$ et $\cL'$  and an exact sequence $0
\to
\cL' \to \cL \oplus \cF \to
\cF' \to 0$, from wich we deduce isomorphisms  $\cExt^i_{\cO_\bP}
(\cF',\cO_\bP)\to\cExt^i_{\cO_\bP}(\cF,\cO_\bP)$ for
$i\geq 2$. The same holds if there is a chain of {\it psi} between $\cF$
and  $\cF'$.

\remk {2.3}  Let $\cF$ and  $\cF'$ two coherent sheaves which are
pseudo-isomorphic. Then the (local) projective dimension of $\cF$ is $\leq
1$ if and only if the same is true for $\cF'$. Therefore we will restrict
ourselves to the set   of  coherent sheaves of 
projective dimension $\leq 1$ , and 
 we will
denote by
${\cal P}\!si$ the  quotient of this set  by the  psi-equivalence relation.
For such a  sheaf
$\cN$,  $H^0\cN(n)$ vanishes for $n \ll 0$, and the first condition of 1.2
is always satisfied.
\vskip 0.3cm

With this notation we have the following result (se also [HMDP]2.4 and
2.8) :

\th {Proposition 2.4}. The injection of the set of vector bundles into the  
set  of  coherent sheaves of
projective dimension $\leq 1$ induces a
natural map from ${\cal S}tab$ to 
${\cal P}\!si$ which is bijective.

\pf
 Suppose that $\cF$ and $ \cF'$ are two pseudo-isomorphic vector bundles such
that  $H_*^2({\cal F})=H_*^2({\cal F}')=0$. By a
``Verdier-lemma'' ([HMDP] 2.11), there exist a coherent sheaf $\cF''$ and two
psi $f:\cF'' \to \cF$ and $f':\cF'' \to \cF'$. By lemma 2.1, we conclude
that 
$\cF''$ is locally free and that $\cF''$ and  $\cF$ (resp. $\cF''$ and 
$\cF'$) are stably isomorphic. Therefore the natural map from ${\cal S}tab$
to ${\cal P}\!si$ is injective.

Now let $\cN$ be a coherent sheaf of
projective dimension $\leq 1$. The graded $R$-modules $N=H_*^0({\cal N})$
and
$\Ext^1_R(N,R)$ are of finite type. Let $L$ be a free
$R$-module of finite type, such that there exists a surjective homomorphism
: $L^{\vee}\to
\Ext^1_R(N,R)$. The composed map :
$$R\to L^{\vee}\otimes L \to \Ext^1_R(N,R)\otimes L \to
\Ext^1_R(N, L)$$
gives an extension of $N$ by $ L$, that is an exact sequence : $0 \to
 L\to F \to N \to 0$. If we dualize, we obtain a coboundary
homomorphism : $L^{\vee}\to \Ext^1_R(N,R)$, which is  the
original surjective homomorphism, so that we get $\Ext^1_R(F, R)=0$. Denoting
by
$\cF$ (resp. $\cL$) the  sheaf associated to $F$ (resp. $L$), we have also
$\Ext^1_{{\cal O}_{\bf P}}(\cF, {\cal O}_{\bf P}(n))=0$ for all $n$. Using
the spectral sequence of $\Ext$ we have also $\cExt^1_{{\cal O}_{\bf P}}(\cF,
{\cal O}_{\bf P})=0$

From the exact sequence : $0 \to
\cL\to \cF \to \cN \to 0$ we deduce isomorphisms  $\cExt^i_{{\cal O}_{\bf
P}}(\cF, {\cal O}_{\bf P})\simeq  \cExt^i_{{\cal O}_{\bf P}}(\cN, {\cal
O}_{\bf P})$ for $i\geq 2$. Since the projective dimension of
$\cN$ is
$\leq 1$, we have $\cExt^i_{{\cal O}_{\bf P}}(\cN, {\cal O}_{\bf
P})=0$ for $i\geq 2$. Hence $\cF$ is a vector bundle, and by the same
argument as in the proof of 2.1, we have also $H^2_*\cF=0$. Therefore the
 map from ${\cal S}tab$ to ${\cal
P}\!si$ is also surjective.

\vskip 0.3cm

Proposition 2.4 proves that there are vector bundles in every class
of ${\cal P}\!si$ and that the psi-equivalence  is an extension to
the set  of  coherent sheaves of (local)
projective dimension $\leq 1$ of the stable equivalence for vector bundles
with $H^2_*=0$.

The following result is an easy consequence of [Ho]. We recall the proof for
the convenience of the reader.

\th {Proposition 2.5}. The application which maps a finite length graded
$R$-module to the sheaf associated to its second module of
syzygies induces a bijective correspondence between
${\cM}_f$ and
${\cal S}tab$.  The application which maps a vector bundle $\cF$ to its first
cohomology module $H^1_*(\cF)$ induces the inverse correspondence between
${\cal S}tab$ and
${\cM}_f$. 

\pf Let $M$ be  a finite length graded
$R$-module, and $F$ be its second module of
syzygies, so that we have an exact sequence :
$0\to F \to L_1 \to L_0 \to M \to 0$
where $L_1$ and $L_0$ are free $R$-modules,
and an exact sequence of (locally free) sheaves :
$0\to \cF \to \cL_1 \to \cL_0 \to 0$.
Taking the associated cohomology exact sequence, we get $M\simeq H^1_*(\cF)$
and $H^2_*(\cF)=0$.

Conversely, let $\cF$ be a vector bundle such that $H^2_*(\cF)=0$, and  let
$L'_1$ be a free
$R$-module of finite type, such that there exists a surjective homomorphism
: $L_1^{\vee}\to H^0_*(\cF^{\vee})$, whose kernel is denoted by $L'_0$. Let
$\cL'_0$ be the associated sheaf. By construction we have $
H^1_*(\cL'_0)=0$, and by Serre duality $ H^1_*(\cF^{\vee})=0$, so we get
also $ H^2_*(\cL'_0)=0$. Therefore by Horrocks theorem $\cL'_0$ is dissoci\'e
; if we set $\cL_0^{'\vee}=\cL_0$, we have an exact sequence : $0\to \cF \to
\cL_1
\to \cL_0 \to 0$ and $\cF$ is the sheaf associated to the second module of
syzygies of $H^1_*(\cF)$.

\vskip 0.3cm
As an immediate consequence of 2.5  we obtain the following result :
  
\th {Corollary 2.6}. The application which maps a  coherent sheaf $\cN$  to
its first cohomology module $H^1_*(\cN)$ induces a bijective correspondence
between ${\cal P}\!si$ and
${\cM}_f$. 

\th {Notation 2.7}. Tensorising a sheaf by ${\cal O}_{\bf P}(n)$ defines an
action of $\bZ$ on
${\cal P}\!si$. We denote by ${\cal P}\!si'$ the quotient  set.

\th {Corollary 2.8}.  The application which maps a curve $C$ to its sheaf
of ideals ${\cal J}_C$ induces a bijective correspondence between ${\cal
B}il$ and
${\cal P}\!si'$.

\pf By the  Rao theorem ([R]), there is a one-to-one
correspondence between the set of curves, up to
biliaison equivalence, and the set of finite length graded $R$-modules,
modulo isomorphism, up to shift in degrees. This correspondence maps the
class of a curve $C$ to the class of its Rao-module $H^1_*({\cal J}_C)$. We
conclude by using 2.6.

\titre {3. Minimal elements}

In the previous paragraph we proved that there are vector bundles in every class
of ${\cal P}\!si$. We will see now that there are
other outstanding elements, namely the reflexive sheaves of  rank 2.
Moreover, we can characterize them.

Recall that for a function $f:\bZ\to\bN$ which vanishes for $m\ll 0$, we
denote by
$f^{\sharp}$ the ``primitive'' function defined by $f^{\sharp}(n)=\sum
_{m\leq n} f(n)$ ([MDP1]I.1).

\th {Proposition 3.1}. ([MDP2]) Let $\cF$ be a vector bundle of rank $r$.
There exists a function $q'_{\cF}:\bZ\to\bN$ with finite support such that
the following properties are equivalent for a dissoci\'e sheaf $\cL=\bigoplus
{\cal O}_{\bf P}(-n)^{l(n)}$ of rank
$r-2$ :\lign
i) a general homomorphism $u:\cL \to \cF$ is injective, and its cokernel is
reflexive, \lign
ii) for all $n$, $l^{\sharp}(n)\leq
(q'_{\cF})^{\sharp}(n)$ , and the condition of the ``obligatory direct
summand'' holds (see [MDP2] for an exact definition).

\remks {3.2} 

1) The result is still true if we suppose that  $\cF$ is only reflexive. In
[MDP2], we give an explicit construction of the function
$q'_{\cF}$. In particular, we have $q'_{\cF\oplus \cL'}=q'_{\cF}+l'$ for
  a dissoci\'e sheaf $\cL'=\bigoplus
{\cal O}_{\bf P}(-n)^{l'(n)}$.

2) It is easy to see that the conditions of 3.1 are satisfied in the case
when $l=q'_{\cF}$. 

\th {Proposition 3.3}.  There exist rank
2 reflexive
 sheaves in every class of ${\cal P}\!si$. The set of first Chern classes of
these reflexive sheaves has a lower bound. All the  rank 2 reflexive sheaves
of the class with minimal $c_1$  (we say that they are {\bf minimal}  rank 2
reflexive sheaves) have the same cohomology.

\pf Let $\cE$ be a  rank 2 reflexive sheaf in the class. In the
last part of the proof of  2.4, we saw that there exists an exact sequence 
: $0 \to
\cL \to \cF \to
\cE \to 0$, where $\cF$ is a vector bundle such that $H^2_*(\cF)=0$ and
where $\cL$ is dissoci\'e.

Let $\cF$ be a vector bundle of rank $r$, $\cL=\bigoplus
{\cal O}_{\bf P}(-n)^{l(n)}$ be a a dissoci\'e sheaf of rank
$r-2$ and $0\to \cL \to \cF \to \cE \to 0$ be an exact sequence where
$\cE$ is reflexive. Let us set  $q'_{\cF}=q'$. We have
$c_1(\cE)=c_1(\cF)+\sum_{n\in {\bf Z} }nl(n)$ and
$\sum_{n\in {\bf Z}}q'(n)=\sum_{n\in
{\bf Z}}l(n)=r-2.$

For
$m\gg 0$ we have : 
$$\sum_{n\in {\bf Z}}nq'(n)=(m+1)\sum_{n\in
{\bf Z}}q'(n)-q^{'\sharp\sharp}(m)$$
$$\sum_{n\in {\bf Z}}nl(n)=(m+1)\sum_{n\in
{\bf Z}}l(n)-l^{\sharp\sharp}(m).$$
By 3.1,  the function $q^{'\sharp}-l^{\sharp}$ is positive, so the function 
$q^{'\sharp\sharp}-l^{\sharp\sharp}$ is positive and increasing. Since it is
positive, we obtain the inequality $\sum_{n\in {\bf Z}}nq'(n)\leq
\sum_{n\in
{\bf Z}}nl(n)$ and $c_1(\cE)\geq c_1(\cF)+\sum_{n\in {\bf Z}}nq'(n)$.

Since it is increasing and vanishes for $m\ll 0$, if it vanishes also for
$m\gg 0$, it vanishes everywhere, and then $l=q'$.

So the set of first Chern classes of rank 2 reflexive sheaves which are
a quotient of
$\cF$ by a dissoci\'e sheaf
has a lower bound, which corresponds to the case when $l=q'$. We say that the
corresponding quotients are {\it minimal rank 2 reflexive quotients} of
$\cF$. 

From the equality $q'_{\cF\oplus \cL'}=q'_{\cF}+l'$ for
a dissoci\'e sheaf $\cL'$, we deduce that the minimal rank 2 reflexive
quotients of $\cF$ are also the minimal rank 2 reflexive
quotients of $\cF\oplus \cL'$. 
Since all the vector bundles of  a given class of ${\cal P}\!si$ are stably
equivalent, they all have the same minimal rank 2 reflexive
quotients. The dimensions of their cohomology groups are all the same, as
follows from the lemma :

\th {Lemma 3.4}. Let $\cF$ be a vector bundle of rank $r$, $\cP$ be a
dissoci\'e sheaf of rank $r-2$, $u$ and $u'$ be two injective homomorphisms 
$\cP\to \cF$ whose cokernels $\cE$ and $\cE'$ are reflexive. Then $\cE$ and
$\cE'$ have the same cohomology, that is $h^i\cE(n)=h^i\cE'(n)$ for all
$i$ and all $ n$.

\pf For all $n$, we have $h^0\cE(n)=h^0\cF(n)-h^0\cP(n)=h^0\cE'(n)$ and 
$h^1\cE(n)=h^1\cF(n))=h^1\cE'(n)$. Moreover, $\cE$ and $\cE'$ have the same
first Chern class $c_1$ and $\cE(n)$ and $\cE'(n)$ have the same Euler
characteristic. Then by Serre duality, $ H^3\cE(n)$ is the dual of ${\rm
Hom\,} (\cE, {\cal O}_{\bf
P}(-n-4)=H^0(\cE^{\vee}(-n-4))=H^0(\cE(-c_1-n-4))$. Therefore for all $n$
$h^3\cE(n)=h^3\cE'(n)$, and since their Euler
characteristic are equal, $h^2\cE(n)=h^2\cE'(n)$.

\remks {3.5}

1) By 2.6 the notion of biliaison of rank 2 reflexive sheaves introduced by
[B] is the restriction of the psi-equivalence to the set of rank 2 reflexive
sheaves. 

2) [B] proves that the minimal  rank 2 reflexive sheaves in a class of
psi-equivalence have also the minimal third Chern class. Hence she points
out that if there are rank 2 vector bundles in the class ($c_3=0$), they are
minimal elements and all the minimal elements are vector bundles.

\vskip 0.3cm

 Now  if we want to give an explicit description
of the inverse correspondence between
${\cal P}\!si'$ and ${\cal B}il$ (cf. 2.8), there are two possibilities,
using the two kinds of outstanding elements in a class of psi-equivalence.

The first choice (a vector bundle) will lead to the theory of
minimal curves ([Mi], [MDP1]~V, [BBM]) and give a characterization of  the
curves in the biliaison class. We use the following result, which is
analogous to 3.1.

\th {Proposition 3.6}. ([MDP2]) Let $\cF$ be a vector bundle of rank $r$.
There exists a function $q_{\cF}:\bZ\to\bN$ with finite support such that
the following properties are equivalent for a dissoci\'e sheaf $\cL=\bigoplus
{\cal O}_{\bf P}(-n)^{l(n)}$ of rank
$r-1$ :\lign
i) a general homomorphism $u:\cL \to \cF$ is injective, and its cokernel is
the twisted ideal of a curve, \lign
ii) for all $n$, $l^{\sharp}(n)\leq
q_{\cF}^{\sharp}(n)$ , and the condition of the ``obligatory direct summand''
holds  (see [MDP2] for
an exact definition).

\remks {3.7} 

1) The function $q_{\cF}$ satisfies also 3.2.1 and 3.2.2, as the function
$q'_{\cF}$.

2) For every curve $C$, there exists  an exact sequence 
: $0 \to
\cL \to \cF \to
\cJ_C \to 0$, where $\cF$ is a vector bundle such that $H^2_*(\cF)=0$ and
where
$\cL$ is dissoci\'e. One can obtain this result as  a consequence of the 
last part of the proof of  2.1, but it has  been already  proved in [MDP1]II
({\it N} resolution of the curve).

3) As in 3.3, we can prove that the curves corresponding to the case when 
$l=q_{\cF}$ have the minimal shift of their ideal, or of their
Rao-module. These are the {\it minimal curves} of the class.

\vskip 0.3cm
With the second choice (a rank 2 reflexive sheaf $\cE$) we obtain curves
which are the zero locus of a section of $\cE(n)$ :
$$0\to {\cal O}_{\bf P}(-n)\to \cE\to {\cJ}_C(c_1(\cE)+n)\to 0.$$

It is then very natural to ask the following question :

\th{Question I}. Is a minimal curve of a biliaison class the zero locus of a
section of a twist of a (minimal in the corresponding class of ${\cal
P}\!si$) rank 2 reflexive sheaf ?

\th {Proposition 3.8}. ([B]) The following properties are equivalent : \lign
i) the general minimal curve  is the zero locus of a
section of a twist of a rank 2 reflexive sheaf,\lign
ii) the general minimal curve is the zero locus of a
section of a twist of a {\bf minimal} rank 2 reflexive
sheaf,\lign
iii) let $\cE$ be a minimal rank 2 reflexive
sheaf and let $n_0$ its minimal twist which has non-zero sections, i.e. $n_0
=
\inf\{\, n\mid H^0\cE(n)\neq 0\,\}$, then  a non-zero section of $\cE(n_0)$
vanishes along a minimal curve,\lign 
iv) for all vector bundle $\cF$ in the
corresponding class, for all $n\in
\bZ$, $q'_{\cF}(n)\leq q_{\cF}(n)$.

It turns out that in general the answer to question I is negative. One
example was given in [B]. We will give a particular case of this example,
where it is possible to give a description of a corresponding
minimal curve.

\exp {3.9} Let $f_1,f_2,f'_3,f'_4$ (resp. $f_1,f_2,f''_3,f''_4$) a regular
sequence of homogeneous polynomials of degrees $n_i=deg(f_i)$ with
$n_1=n_2< n_3=n_4$ and consider the finite
length graded $R$-modules $M'=R/(f_1,f_2,f'_3,f'_4)$,
$M''=R/(f_1,f_2,f''_3,f''_4)$ and $M=M'\oplus M''$. Let $\cF$ be the sheaf 
associated to the second module of
syzygies of $M$ so that we have an exact sequence :
$$0\to \cF \to {\cal O}_{\bf P}(-n_1)^4\oplus {\cal O}_{\bf P}(-n_3)^4 \to
{\cal O}_{\bf P}^2 \to 0.$$ Using the machinery
of [MDP2] we can compute the corresponding functions  
$q_{\cF}$ and $ q'_{\cF}$ and we obtain :
$$q_{\cF}(2n_1)=2,\quad q_{\cF}(n_1+n_3)=3 \quad \hbox {and}
\quad q_{\cF}(n)=0
\quad\hbox {elsewhere},$$
$$ q'_{\cF}(n_1+n_3)=4 \quad \hbox {and}
\quad q'_{\cF}(n)=0
\quad\hbox {elsewhere},$$
which shows that the condition  3.8.iv is not satisfied for $n=n_1+n_3$.

Moreover, any minimal curve $\Gamma$  of the biliaison
class corresponding to $M$ has a resolution :
$$0\to {\cal O}_{\bf P}(-2n_1)^2\oplus {\cal O}_{\bf P}(-n_1-n_3)^3\to \cF
\to \cJ_{\Gamma} (3n_1-n_3)\to 0$$
and a Rao-module $H^1_*(\cJ_{\Gamma})\simeq  M(n_3-3n_1)$.

Let $C'$ (resp. $C''$) be a general minimal curve of the biliaison
class corresponding to $M'$ (resp. $M''$), which has degree $2n_1^2$. The
saturated homogeneous ideal of
$C'$ and $C''$ are given by : $$I_{C'}=(f_1,f_2)^2+(g'_3f'_3+g'_4f'_4) \quad 
I_{C''}=(f_1,f_2)^2+(g''_3f''_3+g''_4f''_4),$$ where $g'_3,g'_4$
(resp. $g''_3,g''_4$) are two independant linear forms in $f_1,f_2$. We have
$H^1_*(\cJ_{C'})=M'(n_3-n_1)$ and
$H^1_*(\cJ_{C''})=M''(n_3-n_1)$.

Now we will use the ``liaison addition'' of Schwartau ([S]). Let $P'$ and
$P''$ be two homogeneous elements of degree
$2n_1$ of the ideal $(f_1,f_2)^2$, without common divisor. The ideal $I=P''
I_{C'}+P'I_{C''}$ is saturated and defines a curve $C$, which is, as a set,
the union of $C'$, $C''$ and the complete intersection defined by $(P',P'')$.
From the exact sequence :
 $$0\to R(-4n_1)\mapright{i} I_{C'}(-2n_1)\oplus
I_{C''}(-2n_1) \mapright{p} I \to 0$$
where the maps $i$ and $p$ are given by  :
$i(a) =(aP', -aP'')$,
$p(b',b'') = b'P''+b''P'$, and the associated exact sequence of sheaves, we
deduce that $C$ has degree
$8n_1^2$ and that its Rao-module is isomorphic to :
$$H^1_*(\cJ_{C'}(-2n_1))\oplus H^1_*(\cJ_{C''}(-2n_1))\simeq
M'(n_3-3n_1)\oplus  M''(n_3-3n_1)\simeq M(n_3-3n_1).$$
Therefore $C$ is a minimal curve  of the biliaison
class corresponding to $M$, and it is not  the zero locus of a
section of a twist of a rank 2 reflexive sheaf. In fact, by 3.8, the general
minimal curve  is not the zero locus of a
section of a twist of a rank 2 reflexive sheaf ; this is an open property,
so no minimal curve  can be the zero locus of a
section of a twist of a rank 2 reflexive sheaf.

\vskip 0.3cm 

By 3.5.2,  if there are rank 2 vector bundles in a class, they
are the minimal elements (among the rank 2 reflexive
sheaves). Recall that   the curves obtained as
the zero locus of a section of a twist of a  rank 2 vector bundle are
exactly the so-called {\it subcanonical} curves (cf. [H]) :

\th {D\'efinition 3.10}. A curve $C$ is subcanonical if there exists $a\in
\bZ$ and an isomorphism : $\omega_C\simeq {\cal O}_C(a)$, where  $\omega_C$
is the dualizing sheaf of $C$.

 In this particular case, the question takes the following form, under which
it was  asked to me by R. Hartshorne and P. Ellia :

\th{Question II}. If  a biliaison class contains subcanonical curves, is a 
minimal curve subcanonical ?

There is a result analogous to 3.8 :

\th {Proposition 3.11}. ([B]) The following properties are equivalent :
\lign i) if  a biliaison class contains subcanonical curves, every 
minimal curve is subcanonical,\lign
ii)  let $\cE$ be a  rank 2 vector bundle  and let $n_0$ its minimal twist
which has non-zero sections, i.e. $n_0 =
\inf\{\, n\mid H^0\cE(n)\neq 0\,\}$, then a non-zero section of $\cE(n_0)$
vanishes along a minimal curve.

In a recent work ([MD]) I gave a positive  answer to question II by proving
3.11.ii.

\vskip 1cm
\titre {References}

[BBM] Ballico E., Bolondi G. Migliore J., The Lazarsfeld-Rao problem for
liaison classes of two-codimensional subschemes of $ {\bf P}^n$. Amer. J.
Math. 113,  117--128 (1991). 

[B] Buraggina A., Biliaison classes of reflexive
sheaves, Math. Nachr. 201, 53--76, 1999.

[H] Hartshorne R., Stable Vector Bundles of rank 2 on $ {\bf P}^3$, Math.
Ann. 238, 229-280 (1978).

[Ho] Horrocks G., Vector bundles on the punctured spectrum of a local
ring, Proc. Lond. Math. Soc., 14, 689-713 (1964).

[HMDP]  Hartshorne R., Martin-Deschamps M. et  Perrin D., Un th\'eor\`eme
de Rao pour les familles de courbes gauches, Journal of Pure and Applied
Algebra  155, 53-76 (2001).

[MD] Martin-Deschamps M., Minimalit\'e des courbes sous-canoniques,
preprint, 2001.

[MDP1]  Martin-Deschamps M. et 
Perrin D., Sur la classification des courbes gauches I, Ast\'erisque, Vol.
184-185, 1990.

[MDP2]  Martin-Deschamps M. et 
Perrin D.,  Quand un morphisme de fibr\'es d\'eg\'en\`ere-t-il le long d'une courbe
lisse ?
		Lecture Notes in Pure and Applied Mathematics Series/200. Marcel Dekker, Inc.
july 1998.

[Mi] Migliore J., Geometric Invariants of Liaison, J.  Algebra 99, 548-572
(1986).

[R] Rao A. P., Liaison among curves in $ {\bf P}^3$, Invent. Math., Vol.
50, 205-217 (1979).

[S] P.Schwartau. Liaison addition and monomial ideals, Thesis, Brandeis
University, Waltham, MA, 1982.

\end